\documentstyle[12pt,amssymb,draft]{amsart}
\theoremstyle{plain}
\newtheorem{thm}{Theorem}[section]
\newtheorem{theorem}[thm]{Theorem}
\newtheorem{corollary}[thm]{Corollary}
\newtheorem{lemma}[thm]{Lemma}
\newtheorem{remark}[thm]{Remark}
\newtheorem{example}[thm]{Example}

\newenvironment{theorem*}[1]{\smallskip\noindent{\bf #1.}\it}{\medskip}
\newenvironment{proofof}[1]{\smallskip\noindent{\it #1}\rm}
                {\hspace*{\fill} $\Box$\medskip}
\newenvironment{proof}{\smallskip\noindent{\it Proof.}\rm}
                        {\hspace*{\fill} $\Box$\medskip}

\numberwithin{equation}{section}
\newcommand{\intd}{\int_{[0,2\pi)}^{\oplus}}
\newcommand\dd{\frac{d^2}{dt^2}}

\newcommand\ov{\overline}

\newcommand\wt{\widetilde}
\newcommand\wh{\widehat}
\newcommand\lan{\langle}
\newcommand\ran{\rangle}
\renewcommand\Re{\operatorname{Re}}

\newcommand\bn{\binom}
\newcommand\tr{\operatorname{tr}}
\newcommand\supp{\operatorname{supp}}
\newcommand\const{\operatorname{const}}

\newcommand\al{\alpha}
\newcommand\be{\beta}
\newcommand\ga{\gamma}

\newcommand\de{\delta}
\newcommand\De{\Delta}
\newcommand\eps{\varepsilon}
\newcommand\la{\lambda}
\newcommand\La{\Lambda}

\newcommand\si{\sigma}
\newcommand\th{\theta}

\newcommand\cH{{\mathcal H}}

\newcommand\bC{{\mathbb C}}

\newcommand\bN{{\mathbb N}}
\newcommand\bR{{\mathbb R}}
\newcommand\bZ{{\mathbb Z}}

\newcommand\fD{{\frak D}}

\newcommand\op{operator}
\newcommand\sa{selfadjoint}

\newcommand\nbh{neighbourhood}

\title{Schr\"odinger operators with periodic singular potentials${}^{\dag}$}
\thanks{${}^{\dag}$The work was partially supported by Ukrainian State Foundation
for Basic Research}

\author{R.~O.~Hryniv and Ya.~V.~Mykytyuk}

\address{Institute for Applied Problems of Mechanics and Mathematics,
	3b~Naukova str., 79601 Lviv, Ukraine}
\email{hryniv@@mebm.lviv.ua}

\address{Department of Mechanics and Mathematics, Lviv National University,
	1 Universytetska str., 79602 Lviv, Ukraine}
\email{mykytyuk@@email.lviv.ua}

\date{\today}

\begin{document}

\begin{abstract}
We show that formal Schr\"odinger operators with singular potentials from
the space $W^{-1}_{2,unif}(\bR)$ can be naturally defined to give
selfadjoint and bounded below operators, which depend continuously
in the uniform resolvent sense on the potential in
the $W^{-1}_{2,unif}(\bR)$-norm.
In the case of periodic singular potentials we also establish pure absolute
continuity and a band and gap structure of the spectrum thus generalising
some classical results for singular potentials of one-dimensional
quasicrystal theory.
\end{abstract}

\maketitle

\section{Introduction}\label{sec:intr}

In a series of recent papers~\cite{SS}--\cite{HM} Schr\"odinger operators
have been considered with singular potentials that are distributions from
the space $W^{-1}_{2,loc}(\bR)$.
More exactly, for a potential $q = \si' + \tau \in W^{-1}_{2,loc}(\bR)$
with real-valued $\si \in L_{2,loc}(\bR)$ and $\tau \in L_{1,loc}(\bR)$
the corresponding Schr\"odinger operator
\begin{equation}\label{eq:S}
	S = - \dd + q
\end{equation}
is defined through
\begin{equation}\label{eq:Sact}
	S u = l (u) := - (u' -\si u)' - \si u' + \tau u
\end{equation}
on the domain
\begin{equation}\label{eq:Sdom}
	\fD(S) = \{ u \in W^1_{1,loc}(\bR) \mid
		u^{[1]}:= u'-\si u \in W^1_{1,loc}(\bR),\
			l(u) \in L_2(\bR)\}.
\end{equation}
It is easily seen that $l(u) = -u'' + qu$ in the sense of
distributions, which implies, firstly, that the operator~$S$ does not
depend on the particular choice of $\si \in L_{2,loc}(\bR)$ and
$\tau\in L_{1,loc}(\bR)$ in the decomposition $q=\si'+\tau$ and,
secondly, that for regular potentials $q \in L_{1,loc}(\bR)$ the
above definition coincides with the classical one. Moreover,
the \op~$S$ is shown to be selfadjoint and bounded below if
$\si$ is compactly supported and $\tau$ is in the limit point case
at~$\pm\infty$~\cite{SS} or if $q \in W^{-1}_2(\bR)$~\cite{NS}.

The regularisation by quasi-derivatives procedure was first suggested
in~\cite{AEZ} for the potential $1/x$ on a finite interval
(see also the books~\cite{N} and \cite{EM} for a detailed exposition
of quasi-differential operators). The more general
setting~\eqref{eq:S}--\eqref{eq:Sdom} developed in~\cite{SS} allows
to consider, e.g., very important cases of Coulomb $1/x$-
and Dirac $\de$-like potentials used to model short- and zero-range
interactions in quantum mechanics
by taking $\si(x) = \log|x|$ and $\si(x) = \chi(x)$,
the Heaviside function, respectively. These two models as well as their
generalizations to potentials that are {\em singular} (i.e., not locally
integrable) on a discrete set were treated in many works, see, e.g.,
\cite{AEZ}, \cite{GK}--\cite{BDL} and the references therein.
In the present paper the potential $q \in W^{-1}_{2,loc}(\bR)$ is not
assumed locally integrable anywhere, though the singularities cannot
be too strong (say, $\de'$-interactions are not allowed).

We remark that it has been realised for a long time that differential
expressions of~\eqref{eq:S} with singular potentials do not generally
determine a unique operator in $L_2(\bR)$.
However, the operator~$S$ defined by~\eqref{eq:Sact}--\eqref{eq:Sdom}
appears to be a ``natural'' \sa\ \op\ associated with~\eqref{eq:S}
for a potential $q \in W^{-1}_2(\bR)$ in the sense that if $q_n$ is any
sequence of regular (infinitely smooth say) potentials that converges to $q$
in $W^{-1}_2(\bR)$, then the corresponding Schr\"odinger operators $S_n$
converge to $S$ in the uniform resolvent sense.
This fact is established in~\cite{SS} for regular Sturm-Liouville operators
on a finite interval or Schr\"odinger \op s with compactly supported~$\si$
and $\tau$ in the limit point case at~$\pm\infty$
and in~\cite{NS} for $q \in W^{-1}_{2}(\bR)$ or for a more general situation
of polyharmonic operators in $\bR^n$ with potentials from an appropriate
space of multipliers. See also~\cite{BFT} and \cite{ENZ} for convergence
results for Schr\"odinger operators with potentials that are Radon measures
and with $\de'$-potentials, respectively, and \cite{AK} for an abstract
setting of form-bounded singular perturbations.

The main aim of this note is to study Schr\"odinger operators with
singular potentials from the space $W^{-1}_{2,unif}(\bR)$,
and in particular with {\em periodic} singular potentials. While this problem
apparently did not receive much attention in the above-cited works
(potentials from the class $W^{s}_{2,unif}(\bR)$, $s>-1$, were considered
in~\cite{He}), the particular cases of periodic and quasiperiodic
$\de$-interactions were quite well understood within the framework of
quasicrystal theory in quantum mechanics, cf.~Kronig-Penney theory and
its various generalizations in~\cite[Ch.~III.2]{AGHH}.
For instance, the Schr\"odinger
operators with periodic $\de$-interactions are shown to have purely
absolutely continuous band spectrum (see also \cite{GK} and \cite{MS} for
more general cases of periodic singular point-like interactions).
Moreover, if $q_n$ is a sequence of regular
short-range interactions that converges to a sum $q$ of
$\de$-interactions in the sense of quadratic forms (which incidentally
implies $W^{-1}_{2,unif}(\bR)$-convergence), then the corresponding
Schr\"odinger \op s $S_n$ converge to $S$ in the uniform resolvent sense.
These two results are generalized in the present paper to an arbitrary
potential from $W^{-1}_{2,unif}(\bR)$.

The main results of our paper are as follows. In Section~\ref{sec:str}
we define the space $W^{-1}_{2,unif}(\bR)$ and show that any real-valued
$q \in W^{-1}_{2,unif}(\bR)$ can be represented (not uniquely) in the form
$q = \si' + \tau$, where $\si$ and $\tau$ are real-valued functions from
$L_{2,unif}(\bR)$ and  $L_{1,unif}(\bR)$, respectively, i.e.,
\begin{align*}
  \|\si\|^2_{2,unif} &:= \sup_{t\in\bR} \int_t^{t+1} |\si(s)|^2 ds < \infty,\\
   \|\tau\|_{1,unif} &:= \sup_{t\in\bR} \int_t^{t+1} |\tau(s)|\,ds < \infty,
\end{align*}
and the derivative is understood in the sense of distributions.
Then in Section~\ref{sec:sa} we prove that the \op~$S$ with
$q = \si' + \tau\in W^{-1}_{2,unif}(\bR)$
as defined by~\eqref{eq:Sact}--\eqref{eq:Sdom} coincides with the form-sum
operator in~\eqref{eq:S} and hence is \sa\ and bounded below.
In Section~\ref{sec:cont} we establish the uniform resolvent convergence result
for $W^{-1}_{2,unif}(\bR)$-convergence of potentials. Finally,
in the last section we consider
a periodic singular potential $q \in W^{-1}_{2,unif}(\bR)$ and prove that
the corresponding \op~$S$ has an absolutely continuous spectrum.

Our results can be illustrated by the following model example.

\begin{example}
Consider the operator $-\De_{\al,\La}$ of the Kronig-Penney theory with
periodic lattice $\La = \{na \mid n \in\bZ \}$~\cite[Ch.~III.2.3]{AGHH}; it
corresponds formally to the potential
\[
	q = \sum_{n\in\bZ} \al\de(\cdot - na)
\]
and is defined rigorously by acting as $-\dd$ on the domain
\[
	\fD(-\De_{\al,\La}) =
		\{ u \in W^1_2(\bR)\cap W^2_2(\bR\setminus \La) \mid
		 u'(na+)-u'(na-) = \al u(na),\ n\in\bZ \}.
\]
We represent $q$ above as $\si' + \tau$ with $\tau \equiv \al/a$ and
$a$-periodic $\si$ equal to $\al/2-\al t/a$ on $[0,a)$; then the corresponding
operator $S$ is easily seen to be exactly $-\De_{\al,\La}$. Our result
implies absolute continuity of the spectrum of $-\De_{\al,\La}$. Although
this statement is well known in the Kronig-Penney theory, its proof within
this theory heavily uses an explicit form of the resolvent of $-\De_{\al,\La}$,
which would not be possible for more general periodic
$q \in W^{-1}_{2,unif}(\bR)$.
\end{example}

Throughout the paper $W^s_2(\bR)$, $s\in\bR$, will denote
the standard Sobolev space, $\|\cdot\|$ without any subscript will
always stand for the $L_2(\bR)$-norm and $f^{[1]}$ for the
{\em quasi-derivative} $f' - \si f$ of a function $f$.

\section{Structure of the space $W^{-1}_{2,unif}(\bR)$}\label{sec:str}

We recall that $W^{-1}_2(\bR)$ is the dual space to the Sobolev space
$W^1_2(\bR)$, i.e., it consists of those {\em distributions}~\cite{Sh}
that define continuous functionals on $W^1_2(\bR)$.
With $\lan \cdot,\cdot\ran$ denoting the duality, we have
for $f\in W^{-1}_2(\bR)$
\[
	\|f\|_{W^{-1}_2(\bR)} := \sup_{0\ne\psi\in W^1_2(\bR)}
			\frac{|\lan \psi, f \ran|}{\|\psi\|_{W^1_2(\bR)}}.
\]
The local uniform analogue of this space is defined as follows. Put
\begin{equation}\label{eq:phi}
   \phi(t):=  \left\{ \begin{array}{ll}
		2(t+1)^2\quad &\mbox{if\quad $t\in [-1,-1/2)$},\\
		1-2t^2\quad   &\mbox{if\quad $t\in [-1/2,1/2)$},\\
		2(t-1)^2\quad &\mbox{if\quad $t\in [1/2,1]$},\\
	        0\quad        &\mbox{otherwise}, 	  \end{array}
              \right.
\end{equation}
and $\phi_n(t) := \phi(t-n)$ for $n\in\bZ$. We say that $f$ belongs to
$W^{-1}_{2,unif}(\bR)$ if $f\phi_n$ is in $W^{-1}_{2}(\bR)$ for all $n\in\bZ$
and
\[
	\|f\|_{W^{-1}_{2,unif}(\bR)} := \sup_{n\in\bZ}
		\|f \phi_n\|_{W^{-1}_{2}(\bR)} < \infty.
\]
Our main aim of this section is to prove the following structure theorem.

\begin{theorem}\label{thm:str}
For any $f\in W^{-1}_{2,unif}(\bR)$ there exist functions
$\si \in L_{2,unif}(\bR)$ and $\tau\in L_{1,unif}(\bR)$ such that
$f = \si'+\tau$ and
\begin{equation}\label{eq:norms}
	C^{-1} \bigl(\|\si\|_{2,unif} + \|\tau\|_{1,unif} \bigr) \le
			\|f\|_{W^{-1}_{2,unif}(\bR)} \le
	C      \bigl(\|\si\|_{2,unif} + \|\tau\|_{1,unif} \bigr)
\end{equation}
with some constant $C$ independent of $f$.
Moreover, the function $\tau$ can be chosen uniformly bounded.
\end{theorem}

We say that $f\in W^{-1}_2(\bR)$ vanishes on an open set $U$ if
	$\lan f, \psi\ran = 0$
whenever $\psi \in W^1_2(\bR)$ has its support in $U$. The support $\supp f$
of $f$ is the complement of the largest open set on which $f$ vanishes.
It follows that $\lan f, \psi \ran$ only depends on the values of $\psi$
in a \nbh\ of $\supp f$, i.e., $\lan f, \psi_1 \ran = \lan f, \psi_2 \ran$
whenever $\psi_1 = \psi_2$ on some open set containing $\supp f$.
In particular, for $f$ with compact support the number $\lan f,1 \ran$
can be defined as $\lan f, \psi \ran$ for any test function $\psi$ that
is identically one on a \nbh\ of $\supp f$.

The crux of the proof of Theorem~\ref{thm:str} is the following

\begin{lemma}\label{lem:str}
Suppose that $f \in W^{-1}_2(\bR)$, $\supp f \subset [-1,1]$, and that
$\lan f, 1 \ran = 0$. Then there exists a function $\si \in L_2(\bR)$
with $\supp \si \subset [-1,1]$ such that $f = \si'$ and
$\|\si\| \le C\|f\|_{W^{-1}_2(\bR)}$ for some positive constant $C$
independent of $f$.
\end{lemma}

\begin{proof}
Denote by $\psi_0$ any test function with $\supp \psi_0 \subset (-1,1)$
and $\lan 1, \psi_0 \ran = 1$. We define a distribution $\si$ by the identity
\begin{equation}\label{eq:drv}
	\lan \si, \psi \ran = - \lan f, J\psi \ran
\end{equation}
where $\psi$ runs over all test functions and
\[
	(J\psi)(t) := \int_{-\infty}^t
		\bigl(\psi (s) - \lan 1, \psi \ran \psi_0(s+2) \bigr)\,ds.
\]
It is easily seen that $J\psi(t) = 0$ for all $|t|$ sufficiently large and
therefore $J\psi$ is a test function and $\si$ is well defined. Observe that
$J\psi' = \psi$ for any test function $\psi$, hence~\eqref{eq:drv}
yields  $f = \si'$ by definition.

Suppose next that $\supp \psi \cap [-1,1] = \emptyset$;
then $J\psi \equiv \const$ on some \nbh\ of $[-1,1]$, whence
$\lan f, J\psi \ran = 0$ by assumption and $\supp \si \subset [-1,1]$.

Finally we show that $\si \in L_2(\bR)$ and that
		$\|\si\| \le C\|f\|_{W^{-1}_2(\bR)}$
for some constant~$C$ independent of $f$.
To this end it suffices to show that the operator $J$ acts boundedly from
$L_2[-1,1]$ into $W^1_2(\bR)$ as~\eqref{eq:drv} then implies
\[
	\|\si\| \le \|J\| \|f\|_{W^{-1}_2(\bR)}
\]
and we can take $C = \|J\|$.

Suppose that $\psi$ is a test function with support in $[-1,1]$. Then
$\supp J\psi \subset [-3,1]$ and since
\(
	\bigl|\int_{-\infty}^t \psi(s) \, ds\bigr|^2 \le 2 \|\psi\|^2,
\)
we have
\[
	\|J\psi\|^2 \le 4 \|\psi\|^2 + 8 \|\psi_0\|^2 \|\psi\|^2
\]
and
\[
	\|(J\psi)'\|^2 = \|\psi\|^2 + 2\|\psi_0\|^2 \|\psi\|^2.
\]
This shows that the norm of~$J$ as an operator from $L_2[-1,1]$ into
$W^1_2(\bR)$ does not exceed
\(
	 4 (1 + \|\psi_0\|),
\)
and the proof is complete.
\end{proof}

\begin{proofof}{Proof of Theorem~\ref{thm:str}.}
For a given $f \in W^{-1}_{2,unif}(\bR)$ and $n\in\bZ$ we put
\[
	f_n := f\phi_n - a_n \chi_{[n-1/2,n+1/2)},
\]
where $a_n := \lan f\phi_n,1 \ran$ and $\chi_\De$ is the characteristic
function of an interval~$\De$. Then $f_n$ satisfies the assumptions of
Lemma~\ref{lem:str} with the interval $[-1,1]$ replaced by $[n-1,n+1]$,
and therefore for every $n \in \bZ$ there exists a function
$\si_n \in L_2(\bR)$ such that $f_n = \si_n'$ and
$\supp \si_n \subset [n-1,n+1]$. It is easily seen that with
\[
	\si := \sum_{n\in\bZ} \si_n \quad \mbox{and}\quad
	\tau := \sum_{n\in\bZ} a_n \chi_{[n-1/2,n+1/2)}
\]
we have $f = \si' + \tau$, so it remains to show that $\si$ and $\tau$ belong
to $L_{2,unif}(\bR)$ and $L_{1,unif}(\bR)$ respectively and that
inequality~\eqref{eq:norms} holds.

Denote by $\psi_n$ a $W^1_2(\bR)$-function with support in $[n-2,n+2]$
that is identically one on $[n-3/2,n+3/2]$ and linear on $[n-2,n-3/2]$
and $[n+3/2,n+2]$. Then $\|\psi_n\|_{W^1_2(\bR)} \le 3$, whence
\[
	|a_n| = |\lan f\phi_n, \psi_n\ran| \le 3\|f\phi_n\|_{W^{-1}_2(\bR)}
		\le 3 \|f\|_{W^{-1}_{2,unif}(\bR)}
\]
and
\[
	\|f_n\|_{W^{-1}_2(\bR)} \le \|f\phi_n\|_{W^{-1}_2(\bR)}
		+ |a_n|\|\chi_{[n-1/2,n+1/2}\|_{W^{-1}_2(\bR)}
		\le 4 \|f\|_{W^{-1}_{2,unif}(\bR)}.
\]
The above inequalities yield the estimates
\[
	\|\tau\|_{1,unif} \le \sup |\tau|
		\le 3 \|f\|_{W^{-1}_{2,unif}(\bR)}
\]
and
\[
	\|\si\|_{2,unif} \le 8\|J\| \|f\|_{W^{-1}_{2,unif}(\bR)},
\]
which establishes one part of the inequality required. For the second one we
observe that for $\si \in L_{2,unif}(\bR)$ and $\tau\in L_{1,unif}(\bR)$
it holds
\[
	|\lan \si'\phi_n, \psi \ran| = |\lan \si, (\phi_n\psi)' \ran|
		\le  2 \|\si\|_{2,unif} \|(\phi_n\psi)'\|
		\le  6 \|\si\|_{2,unif} \|\psi\|_{W^1_2(\bR)}
\]
and
\[
	|\lan \tau\phi_n, \psi \ran| \le
		2 \|\tau\|_{1,unif}\sup |\psi|
		\le 2 \|\tau\|_{1,unif} \|\psi\|_{W^1_2(\bR)}.
\]
Here we have used the inequality $\sup |\psi| \le \|\psi\|_{W^1_2(\bR)}$,
which follows from the relations
\[
     |\psi(t)|^2 = 2 \int_{-\infty}^t \Re \psi'\ov \psi\le
	 \int_{-\infty}^t (|\psi'|^2 + |\psi|^2) \le \|\psi\|^2_{W^1_2(\bR)}.
\]
Therefore $\si'$ and $\tau$ also belong to the space $W^{-1}_{2,unif}(\bR)$
and, moreover,
\begin{align*}
	\|\si'\|_{W^{-1}_{2,unif}(\bR)} &\le 6 \|\si\|_{2,unif},\\
	\| \tau\|_{W^{-1}_{2,unif}(\bR)} & \le 2 \|\tau\|_{1,unif},
\end{align*}
and the theorem is proved.
\end{proofof}

\begin{remark}\label{rem:per} \rm
We say that a distribution $f$ is $T$-periodic if
\(
	\lan f, \psi (t) \ran = \lan f , \psi (t+T) \ran
\)
for any test function $\psi$. It is easily seen that for a $1$-periodic
potential $f\in W^{-1}_{2,unif}(\bR)$ the above construction gives a
$1$-periodic function $\si$ and a constant function
$\tau \equiv \lan f\phi_0,1 \ran$. If $f$ is $T$-periodic, we first
apply the construction to the $1$-periodic potential $\hat f(t):= f(Tt)$
to write $\hat f = \hat \si' + \hat\tau$
with $1$-periodic $\hat \si$ and $\hat\tau \equiv \lan\hat f\phi_0,1\ran$,
and then after rescaling we get $f = \si' + \tau$ with $T$-periodic
$\si(t) := T \hat\si (t/T)$ and $\tau := \hat \tau$.
\end{remark}

\section{Selfadjointness of the operator~$S$}\label{sec:sa}

\def\t{{\frak t}}
In this section, we shall prove that the operator~$S$ as given
by~\eqref{eq:Sact} and \eqref{eq:Sdom} is selfadjoint and bounded below.
In fact, we shall show that the quadratic form of the operator~$S$
coincides with
\[
	\t(u):= (u',u') - (\si u',u) - (\si u,u') + (\tau u,u)
\]
and the latter is a relatively bounded perturbation of the form
$\t_0(u):=(u',u')+(u,u)$ with relative bound zero.
Therefore $t$ is closed and bounded below on the domain $W^1_2(\bR)$,
whence $S$ is a selfadjoint bounded below \op\ and
$\fD(S)\subset\fD(\t)=W^1_2(\bR)$.

\begin{lemma}\label{lem:ineq}
For any $f \in W^1_2[0,1]$ and any $\eps \in(0,1]$ the following
inequalities hold:
\begin{align}\label{eq:max}
	\max_{t\in[0,1]}|f(t)|^2 \ \
		&\le \eps \int_0^1 |f'|^2 dt + 8\eps^{-1} \int_0^1|f|^2dt,\\
	\Bigl(\int_0^1 |f'\ov f|^2 dt\Bigr)^{1/2}
		&\le \eps \int_0^1 |f'|^2 dt + 4\eps^{-3} \int_0^1|f|^2dt.
	\label{eq:f'f}
\end{align}
\end{lemma}

\begin{proof}
For an arbitrary function $\phi \in W^1_2(\bR)$ and any $\eta>0$ we find that
\[
	|\phi(t)|^2 = \int_{-\infty}^t\frac{d}{ds}|\phi(s)|^2ds
		=   2\int_{-\infty}^t \Re \phi'\ov\phi\,ds
		\le \eta \|\phi'\|^2 + \eta^{-1}\|\phi\|^2.
\]
Given a function $f \in W^1_2[0,1]$, we extend it to $\phi\in W^1_2(\bR)$
through
\[
	\phi(t) := \left\{ \begin{array}{ll}
		f(t)\quad   	   &\mbox{if\quad $t\in [0,1]$},\\
		f(2-t)(2-t)\quad   &\mbox{if\quad $t\in (1,2]$},\\
		f(-t)(1+t)\quad    &\mbox{if\quad $t\in [-1,0)$},\\
	        0\quad             &\mbox{otherwise};
		           \end{array}    \right.
\]
then
\[
	\|\phi' \|^2 \le 3 \|f'\|^2_{L_2(0,1)} + 4 \|f\|^2_{L_2(0,1)}
		\quad\text{and}\quad
	\|\phi\|^2 \le 2 \|f\|^2_{L_2(0,1)}.
\]
Therefore
\[
	\max_{t\in[0,1]}|f(t)|^2 \le \eta\|\phi'\|^2 + \eta^{-1}\|\phi\|^2
		\le 3\eta \|f'\|^2_{L_2(0,1)}
		+(4\eta + 2\eta^{-1}) \|f\|^2_{L_2(0,1)},
\]
which implies~\eqref{eq:max} upon setting $\eps = 3\eta \le 1$.
Using~\eqref{eq:max} with $\eps^2$ instead of $\eps$,
we derive the inequality
\begin{align*}
	\|f'\ov f\|^2_{L_2(0,1)} \le \max_{t\in[0,1]}|f(t)|^2\|f'\|^2_{L_2(0,1)}
		&\le \eps^2 \|f'\|^4_{L_2(0,1)}
			+ 8\eps^{-2}\|f'\|^2_{L_2(0,1)}\|f\|^2_{L_2(0,1)}\\
		&\le \bigl(\eps \|f'\|^2_{L_2(0,1)}
			+ 4\eps^{-3}\|f\|^2_{L_2(0,1)}\bigr)^{2},
\end{align*}
and the proof is complete.
\end{proof}

\begin{lemma}\label{lem:tform}
The quadratic form $\t$ is closed and bounded below on~$W^1_2(\bR)$.
\end{lemma}

\begin{proof}
Suppose that $\si\in L_{2,unif}(\bR)$, $\tau \in L_{1,unif}(\bR)$, and
$u \in W^1_2(\bR)$. Using the relations~\eqref{eq:max} and \eqref{eq:f'f},
we find that, with an arbitrary $\eps\in(0,1]$ and $\eta\in(0,1]$,
\begin{align*}
   |(\si u',u)| \le \sum_{n\in\bZ} \int_n^{n+1}|\si u'\ov u|\,dt
		&\le \sum_{n\in\bZ} \Bigl(\int_n^{n+1}|\si|^2dt\Bigr)^{1/2}
				\Bigl(\int_n^{n+1} |u'u|^2dt\Bigr)^{1/2} \\
		&\le \|\si\|_{2,unif} (\eps\|u'\|^2 + 4\eps^{-3}\|u\|^2)
\end{align*}
and
\begin{align*}
	|(\tau u,u)| \le \sum_{n\in\bZ} \int_n^{n+1}|\tau u \ov u|\,dt
		&\le \sum_{n\in\bZ} \int_n^{n+1}|\tau|\,dt
			\max_{t\in[n,n+1]} |u(t)|^2\\
	&\le \|\tau\|_{1,unif} (\eta\|u'\|^2 + 8\eta^{-1}\|u\|^2).
\end{align*}
This shows that the quadratic form $(\si u',u) + (\si u,u') + (\tau u,u)$
is bounded with respect to the form $\t_0$ with relative bound zero.
Therefore by the KLMN theorem (see~\cite[Theorem~X.17]{RS2}) the quadratic
form~$\t$ is closed and bounded below on the domain
$\fD(\t) = \fD(\t_0) = W^1_2(\bR)$. The lemma is proved.
\end{proof}

\begin{remark}\label{rem:lbnd}\rm
Using the above inequalities with $\eps = \min\{1, (4\|\si\|_{2,unif})^{-1}\}$
and $\eta = \min\{ 1, (2\|\tau\|_{1,unif})^{-1}\}$, we find that the quadratic
form $\t$ is bounded below by
\[
	\ga(\t) := - \bigl( 2(4\|\si\|_{2,unif})^4 +
			      16\|\tau\|^2_{1,unif}+ 6 \bigr).
\]
Recalling~Theorem~\ref{thm:str} we can recast this as
\begin{equation}\label{eq:lbound}
	\ga(\t) \ge - \bigl(a \|q\|_{W^{-1}_{2,unif}(\bR)} + b \bigr)^{4}
\end{equation}
with some $a,b>0$ independent of $q$.
\end{remark}

\medskip

Denote by $T$ a \sa\ operator that is associated with the form~$\t$
according to the second representation theorem (see, e.~g.,
\cite[Theorem~VI.2.23]{Ka}). Recall that $T$ is the \sa\ operator
for which $\fD(T)\subset \fD(\t)$ and the equality
$(Tu,v) = \t(u,v)$ holds for all $u \in \fD(T)$ and all $v\in\fD(\t)$.

\begin{theorem}\label{thm:S=T}
The operator $S$ coincides with $T$. In particular, the operator $S$ is
selfadjoint, bounded below, and $\fD(S) \subset W^1_2(\bR)$.
\end{theorem}

\begin{proof}
Fix $u\in\fD(T)$ and take an arbitrary $v \in \fD(\t)$. Then $u \in W^1_2(\bR)$
and by~\eqref{eq:max} with $\eps =\sqrt8$
\begin{multline*}
 \qquad	\int_\bR |\si u|^2 \le \sum_{n\in\bZ} \int_n^{n+1} |\si u|^2 \le
		\|\si\|^2_{2,unif} \sum_{n\in\bZ} \max_{t\in[n,n+1)} |u(t)|^2
	\\      \le \sqrt8 \|\si\|^2_{2,unif}\bigl(\|u'\|^2 + \|u\|^2\bigr)
		< \infty.\qquad
\end{multline*}
Therefore $\si u \in L_2(\bR)$ and $u^{[1]} \in L_2(\bR)$; now
\begin{multline*}
 \qquad	(Tu,v) = \t(u,v) = (u',v') - (\si u,v') - (\si u',v) + (\tau u,v)\\
		=(u^{[1]},v') - (\si u',v) + (\tau u,v)
		= \bigl(-(u^{[1]})' - \si u' + \tau u,v\bigr) \qquad
\end{multline*}
and hence $T u = -(u^{[1]})' - \si u' + \tau u$ in the sense of distributions.
Observe that $T u \in L_2(\bR) \subset L_{1,loc}(\bR)$ and
$\si u', \tau u\in L_{1,loc}(\bR)$; it follows that
$(u^{[1]})' \in L_{1,loc}(\bR)$ and so $u^{[1]} \in W^1_{1,loc}(\bR)$.
Therefore $u \in \fD(S)$ and $T \subset S$; since $S$ is evidently a
symmetric operator, we conclude that $T=S$. The proof is complete.
\end{proof}

\section{Continuous dependence on the potential}\label{sec:cont}

In this section, we shall prove that the Schr\"odinger operator $S$
defined by~\eqref{eq:Sact}--\eqref{eq:Sdom} depends continuously
(in the sense of the uniform resolvent convergence) on the potential~$q$
in the $W^{-1}_{2,unif}(\bR)$-norm. We remark that this generalizes
the convergence results with respect to the topology of $W^{-1}_{2}(\bR)$
of~\cite{NS} and the $\star$-weak topology of Radon measures of~\cite{BFT}.

\begin{theorem}\label{thm:cnt}
Suppose that $q_n \in W^{-1}_{2,unif}(\bR)$, $n\in\bN$, is a sequence
of potentials that converges to $q$ in $W^{-1}_{2,unif}(\bR)$-norm
and that $S_n$ and $S$ are the corresponding Schr\"odinger operators.
Then $S_n$ converge to $S$ as $n \to \infty$
in the uniform resolvent sense, i.~e.,
\[
	\|(S_n - \la)^{-1} - (S - \la)^{-1}\| \to 0
		\qquad\text{as}\qquad	n \to \infty
\]
for any $\la$ in the resolvent set of $S$ and $S_n$, $n\in\bN$.
\end{theorem}

\begin{proof}
Observe first that the corresponding quadratic forms $\t_n$ and $\t$ have
the same domain $W_2^1(\bR)$ and are uniformly bounded below
(see~\eqref{eq:lbound}). We choose sequences $\si_n \in L_{2,unif}(\bR)$ and
$\tau_n \in L_{1,unif}(\bR)$ so that $q_n -q = \si'_n + \tau_n$ and
\[
	\|\si_n\|_{2,unif} + \|\tau_n\|_{1,unif} \le
		C \|q_n - q\|_{W^{-1}_{2,unif}(\bR)}
\]
with the constant $C$ of Theorem~\ref{thm:str}. Repeating the arguments
of the proof of Lemma~\ref{lem:tform} we find that
\begin{multline*}
	|\t_n(u) - \t(u)| \le |(\si_nu',u)| + |(\si_nu,u')| + |(\tau_nu,u)|\\
		\le C_1\bigl(\|\si_n\|_{2,unif} + \|\tau\|_{1,unif}\bigr)
		\|u\|^2_{W^1_2(\bR)}
		\le CC_1 \|q_n - q\|_{W^{-1}_{2,unif}(\bR)} \|u\|^2_{W^1_2(\bR)}
\end{multline*}
for all $u\in W_2^1(\bR)$ and some constant $C_1>0$. The claim follows now
from~\cite[Theorem~VIII.25c]{RS1}, and the proof is complete.
\end{proof}

\begin{remark}\label{rem:mult} \rm
An alternative way to prove Theorem~\ref{thm:cnt} is to use the results of
Section~\ref{sec:sa} to show that the space $W^{-1}_{2,unif}(\bR)$ is
embedded into the space of multipliers $M_0[1]$; then the claim of
Theorem~\ref{thm:cnt} follows from Theorem~9 of~\cite{NS}.
\end{remark}

\section{Periodic potentials}\label{sec:per}

Suppose now that the potential $q\in W^{-1}_{2,unif}(\bR)$ is $1$-periodic;
then (recall Remark~\ref{rem:per}) $q = \si' + \tau$ with $1$-periodic
$\si\in L_{2,unif}(\bR)$ and $\tau \equiv \lan q \phi_0, 1 \ran$.
The purpose of this section is to show that in this case the spectrum of $S$
is absolutely continuous and has a band and gap structure.

Proof of absolute continuity of~$S$ follows the standard spectral analysis of
periodic Schr\"odinger operators (cf.~\cite[Ch.~XIII.16]{RS4}).
We decompose the space $L_2(\bR)$ into a direct integral
\[
	\intd \cH'\,\frac{d\th}{2\pi} =:\cH
\]
with identical fibres $\cH':= L_2[0,1)$; then the operator
$U: L_2(\bR) \to \cH$ defined by
\begin{equation}\label{eq:U}
  (U f)(x,\th) = \sum_{n=-\infty}^{\infty} e^{-in\th} f(x+n)
\end{equation}
is unitary. Now the \op\ $\wt S:= U S U^{-1}$ is unitarily equivalent
to $S$ and can be decomposed into the direct integral
\begin{equation}\label{eq:dcmp}
  \wt S = \intd S_\th\, \frac{d\th}{2\pi},
\end{equation}
where $S_\th$ is the \op\ in $\cH'$ defined by
\[
	S_\th f = l (f)
\]
on the domain
\begin{multline*} \qquad\quad
	\fD(S_\th) = \{ f \in W^1_{1,loc}[0,1) \mid
		f^{[1]} \in W^1_{1,loc}[0,1),\ l(f) \in \cH' ,\\
		f^{[1]}(1) = e^{i\th}f^{[1]}(0),
		f(1) = e^{i\th} f(0) \}. \qquad\quad
\end{multline*}
To show this, we denote by $\wh S$ the operator given by the right
hand side of~\eqref{eq:dcmp}. For any compactly supported $f \in \fD(S)$
the sum in~\eqref{eq:U} is finite and hence the relations $Uf \in \fD(\wh S)$
and $\wh SUf = USf$ are straightforward. Next we observe that the functions
$f \in \fD(S)$ with compact support constitute a core of $S$.
Therefore for any $f \in \fD(S)$ there exist compactly supported
$f_n\in\fD(S)$ such that $f_n \to f$ and $Sf_n \to Sf$ in $L_2(\bR)$ as
$n\to \infty$; then $Uf_n \to Uf$ and $\wh SUf_n = USf_n\to USf$ in $\cH$
as $n\to\infty$. As $\wh S$ is closed we get $Uf\in\fD(\wh S)$ and
$\wh SUf = USf$; therefore $\wt S = USU^{-1} \subset \wh S$ and
$\wt S = \wh S$ since both these operators are selfadjoint.

It follows from decomposition~\eqref{eq:dcmp} that a number $\la$
belongs to the spectrum $\si(\wt S)$ of the operator $\wt S$ if and only if
for any $\eps > 0$
\begin{equation}\label{eq:spctr}
 d\mu\bigl\{ \th \in [0,\pi]
        \mid \si(S_\th) \cap (\la - \eps, \la + \eps) \ne \emptyset \bigr\} >0,
\end{equation}
where $d\mu$ denotes the {\em Lebesgue measure\/} on $\bR$
(\cite[Theorem~XIII.85d]{RS4}).

Observe~\cite{SS} that all operators $S_\th$ have discrete spectra.
We shall prove the following result.

\begin{lemma}\label{lem:anal}
For every fixed nonreal $\la$ the resolvent $(S_\th - \la)^{-1}$ is an
analytic operator function of~$\th$  in a \nbh\ of $(0,2\pi)$.
\end{lemma}

\begin{proof}
Denote by $u_1 = u_1(t,\la)$ and $u_2=u_2(t,\la)$ solutions of equation
$l(u) = \la u$ satisfying the boundary conditions $u_1(0)=0$ and $u_2(1) =0$.
We recall that $u$ being a solution of $l(u) = \la u$ means that
\begin{equation}\label{eq:syst}
	\frac{d}{dt} \bn{u^{[1]}}{u} = \begin{pmatrix}
		 -\si & -\si^2 + \tau - \la \\ 1 & \si
	\end{pmatrix} \bn{u^{[1]}}{u}
\end{equation}
and hence $u$ enjoys standard uniqueness properties of solutions to
second order differential equations with regular (i.e. locally integrable)
coefficients.

Observe that $u_1(1)\ne 0$ and $u_2(0) \ne 0$ as otherwise $\la$ would be
an eigenvalue of the operator $S_D$ determined in $\cH'$ by the differential
expression $l$ and the Dirichlet boundary conditions~\cite{SS}; since
$S_D$ is a selfadjoint \op, this is impossible. Moreover, it follows from
the Cayley-Hamilton theorem that the Wronskian
\[
	W(t) := u_1(t) u_2^{[1]}(t) -  u_1^{[1]}(t) u_2(t)
\]
does not depend on $t\in [0,1]$ and hence can be normalized to be $1$.
It is easily seen that the resolvent $(S_D-\la)^{-1}$ is given then by an
integral operator $(K_\la f)(t) = \int_0^1 K(t,s) f(s)\,ds$ with the kernel
\[
	K(t,s) = \left\{ \begin{array}{ll}
		u_1(t)\ov u_2(s)\quad &\mbox{if\quad $t\le s$},\\
		\ov u_1(s)u_2(t)\quad &\mbox{if\quad $t\ge s$}.
		\end{array}          \right.
\]

Consider the difference $v:=(S_\th - \la)^{-1}f - (S_D -\la)^{-1}f$;
this function solves the equation $l (v) = \la v$ and hence equals
$\al_1(\th,f)u_1 + \al_2(\th,f)u_2$ for some coefficients $\al_1$ and $\al_2$
dependent on $\th$ and $f$. The function
\[
	w:= (S_\th - \la)^{-1}f = (S_D - \la)^{-1}f
			+ \al_1(\th,f)u_1 + \al_2(\th,f)u_2
\]
satisfies the quasiperiodic boundary conditions $w(1) = e^{i\th}w(0)$
and $w^{[1]}(1) = e^{i\th} w^{[1]}(0)$, which means that $\al_1$ and $\al_2$
must solve the system
\begin{align*}
	\al_1 u_1(1) - e^{i\th} \al_2 u_2(0) &=0;\\
	\al_1\Bigl\{ u^{[1]}_1(1) - e^{i\th}u_1^{[1]}(0)\Bigr\} +
	\al_2\Bigl\{ u^{[1]}_2(1) - e^{i\th}u_2^{[1]}(0)\Bigr\} &=
	\be(f,\th),
\end{align*}
where $\be(f,\th) :=e^{i\th}u_1^{[1]}(0) \int_0^1 f \ov u_2 -
		u_2^{[1]}(1) \int_0^1 f \ov u_1$.
We observe that for $f=0$ and any $\th \in [0,2\pi)$ the homogeneous system
above has only the trivial solution $\al_1 = \al_2 = 0$ as otherwise the
function $w = \al _1u_1 + \al_2 u_2$ would satisfy the quasiperiodic boundary
conditions and hence the nonreal number $\la$ would be an eigenvalue of the
selfadjoint operator $S_\th$, which is impossible. Therefore the
discriminant~$d(\th)$ of the above system is nonzero and its solution equals
\[
	\al_1(f,\th) = \frac{e^{i\th}u_2(0)}{d(\th)}\be(f,\th), \qquad
	\al_2(f,\th) = \frac{u_1(1)}{d(\th)}\be(f,\th).
\]
We see that $\al_1$ and $\al_2$ are continuous linear
functionals of $f$ depending analytically on $\th \in (0,2\pi)$. This
completes the proof of analyticity of $(S_\th - \la)^{-1}$.
\end{proof}

It follows from Lemma~\ref{lem:anal} that the eigenvalues $\la_k(\th)$
and the eigenvectors $v_k(\th)$, $k\in \bN$, can be labelled to be analytic
in~$\th$. By~\eqref{eq:spctr} the spectrum of $S$ is now just the union
of ranges of the functions $\la_k(\th)$, $k\in\bN$, when $\th$ varies over
$[0,2\pi)$; this establishes the so-called band and gap structure of $\si(S)$.
Absolute continuity of $S$ will follow
from~\cite[Theorem~XIII.86]{RS4} as soon as we show that all $\la_k(\th)$
are nonconstant. To this end we shall give an alternative description of
the spectra of $S_\th$.

Denote by $v_1(t,\la)$ and $v_2(t,\la)$ solutions of the
system~\eqref{eq:syst} satisfying the following initial conditions:
\begin{align*}
	v_1^{[1]}(0,\la) = 1, \qquad & v_1(0,\la) = 0, \\
	v_2^{[1]}(0,\la) = 0, \qquad & v_2(0,\la) = 1.
\end{align*}
Then the {\em fundamental matrix} $M(t,\la)$ given by
\[
	M(t,\la) : = \begin{pmatrix}
		v_1^{[1]}(t,\la) & v_2^{[1]}(t,\la)\\
		v_1 (t,\la)      & v_2 (t,\la)
	\end{pmatrix}
\]
depends continuously on $t\in[0,1]$ and analytically on $\la \in \bR$,
$\det M(t,\la) \equiv 1$ by the Cayley-Hamilton theorem,
and for any solution $X(t)$ of equation~\eqref{eq:syst} the following
equality holds
\[
	X(t) = M(t,\la)X(0).
\]
If $X(t) = (x^{[1]}(t),x(t))^T$ and $x(t)$ is an eigenfunction of $S_\th$,
then the vector $X$ satisfies the boundary condition $X(1) = e^{i\th}X(0)$,
which implies that $M(1,\la)X(0) = e^{i\th} X(0)$. Therefore $\la \in \bR$
is an eigenvalue of the \op~$S_\th$ if and only if $e^{i\th}$ is an
eigenvalue of the matrix $M(1,\la)$. Since $\det M(1,\la) \equiv 1$, the
latter condition is equivalent to the equality
\begin{equation}\label{eq:trace}
	\tr M(1,\la) = 2 \cos \th.
\end{equation}

\begin{lemma}\label{lem:monot}
The function $\tr M(1,\la)$ is strictly monotone at a point $\la_0$
whenever $|\tr M(1,\la_0)| < 2$.
\end{lemma}

\begin{proof}
The statement of the lemma will follow from the results of~\cite{HSV}
as soon as we show that the matrix $M(1,\la)$ is {\em positively rotating}.
We recall that this requires the function $\arg\Bigl( M(1,\la) X\Bigr)$
to be strictly increasing in $\la \in \bR$ for any nonzero vector
$X \in \bC^2$; here for a vector $X=(x_1,x_2)$ we put
$\arg X := \arg (x_1 + i x_2)$ measured continuously in $X$.

Observe that if
$X(t,\la) = \bigl( x_1(t,\la), x_2(t,\la)\bigr) := M(t,\la) X$,
then by definition $X(t,\la)$ solves the system~\eqref{eq:syst}
and hence $x_2$ satisfies the equation
$l(u) = \la u$ and $x_1(t,\la) = x_2^{[1]}(t,\la)$.
Put $\th(t,\la) := \arg X(t,\la)$; then
$\cot\th (t,\la) \equiv {x_2^{[1]}}/{x_2}$.
After differentiating both sides in $t$ we get
\[
     -\frac{\th'}{\sin^2{\th}} = \frac{-\si x_2'x_2 - x_2^{[1]}x_2'}{x_2^2}
				+\tau - \la = -(\cot\th + \si)^2 + \tau - \la,
\]
or
\[
	\th' = \la \sin^2\th - \tau \sin^2\th + (\cos\th + \si\sin\th)^2.
\]
It follows from~\cite[proof of Theorem~XI.3.1]{Ha} that the function
$\th(1,\la)$ is strictly increasing in $\la$, and therefore $M(1,\la)$
is positively rotating and the claim of the lemma follows.
\end{proof}

Recalling now relation~\eqref{eq:trace}, we derive the following

\begin{corollary}\label{cor:monot}
Every eigenvalue $\la_k(\th)$, $k\in\bN$, when chosen continuous in $\th$,
is analytic and strictly monotone in $\th$ on the intervals $(0,\pi)$
and $(\pi,2\pi)$.
\end{corollary}

Now we combine the results obtained and apply \cite[Theorem~XIII.86]{RS4}
to arrive at the following conclusion.

\begin{theorem}\label{thm:per}
Suppose that the potential $q\in W^{-1}_{2,unif}(\bR)$ is periodic and let
$S$ denote the corresponding Schr\"odinger \op\ of~\eqref{eq:S} constructed
by~\eqref{eq:Sact}--\eqref{eq:Sdom}. Then the spectrum of $S$ is purely
absolutely continuous and has a band and gap structure.
\end{theorem}

{\em Acknowledgements.} The authors thank Prof.~V.~A.~Mikhailets for
enlightening discussions and many useful remarks and comments.


\end{document}